\documentclass[12pt]{article}

\usepackage{amsmath, amsthm, amsfonts, amssymb}
\usepackage[all,cmtip]{xy}
\usepackage[margin=2.5cm]{geometry}
\usepackage{parskip}
\usepackage{graphicx}
\usepackage{eucal}

\numberwithin{equation}{section}

\theoremstyle{plain}
\newtheorem{thm}{Theorem}[section]

\theoremstyle{definition}
\newtheorem{defn}{Definition}[section]

\DeclareMathOperator{\en}{End}

\DeclareMathOperator{\str}{STr}
\DeclareMathOperator{\tr}{Tr}
\DeclareMathOperator{\zhu}{Zhu}

\newcommand{\C}{\mathbb{C}}

\newcommand{\Q}{\mathbb{Q}}
\newcommand{\Z}{\mathbb{Z}}

\newcommand{\CC}{\mathcal{C}}
\newcommand{\HH}{\mathcal{H}}

\newcommand{\al}{\alpha}

\newcommand{\G}{\Gamma}
\newcommand{\D}{\Delta}

\newcommand{\om}{\omega}
\newcommand{\vp}{\varphi}

\newcommand{\slmat}
{\left(\begin{smallmatrix} a & b \\ c & d \end{smallmatrix}\right)}
\newcommand{\ov}{\overline}

\newcommand{\vac}{{\left|0\right>}}

\newcommand{\NS}{\text{NS}}

\newcommand{\squarething}[1]{
  \noindent
  \begin{center}
  \framebox{
    \vbox{
      \vspace{4mm}
      \hbox to 6.3in { {\Large \hfill #1  \hfill} } 
      \vspace{4mm}
    }
  }
  \end{center}
  \vspace*{4mm}
}

\newcommand{\twobytwo}[4]
{\left(\begin{smallmatrix} #1 & #2 \\ #3 & #4 \end{smallmatrix}\right)}

\usepackage{stmaryrd}

\newcommand{\nstw}{\NS^{\text{tw}}}
\newcommand{\nsuntw}{\NS^{\text{untw}}}

\title{Vertex Operator Superalgebras and Odd Trace Functions}
\date{}
\author{Jethro Van Ekeren\footnote{email: jethrovanekeren@gmail.com}\\ \small{\emph{IMPA, Rio de Janeiro, RJ 22460-320, Brasil}}}

\begin{document}

\maketitle

\begin{abstract}
We begin by reviewing Zhu's theorem on modular invariance of trace functions associated to a vertex operator algebra, as well as a generalisation by the author to vertex operator superalgebras. This generalisation involves objects that we call `odd trace functions'. We examine the case of the $N=1$ superconformal algebra. In particular we compute an odd trace function in two different ways, and thereby obtain a new representation theoretic interpretation of a well known classical identity due to Jacobi concerning the Dedekind eta function.
\end{abstract}

\section{Introduction}

One of the most significant theorems in the theory of vertex operator algebras (VOAs) is the modular-invariance theorem of Zhu \cite{Zhu}. The theorem states that under favourable circumstances the graded dimensions of certain modules over a VOA are modular forms for the group $SL_2(\Z)$. The favourable circumstances are that the VOA be rational, $C_2$-cofinite, and be graded by integer conformal weights (we define all terms in Section \ref{definitions} and state Zhu's theorem fully in Section \ref{zhusection} below).

Numerous generalisations of Zhu's theorem have appeared in the literature: to twisted modules over VOAs \cite{DLMorbifold}, to vertex operator superalgebras (VOSAs) and their twisted modules \cite{DZ}, \cite{DZ2} (see also \cite{Jordan}), to intertwining operators for VOAs \cite{Huang}, \cite{M2}, to twisted intertwining operators \cite{Y}, and to non rational VOAs \cite{Mnonrational}.

In \cite{meCMP} the present author relaxed the assumption of integer conformal weights of $V$ to allow arbitrary rational conformal weights. This work was carried out in the setting of twisted modules over a rational $C_2$-cofinite VOSA. Actually it is worth noting that in that paper the condition of $C_2$-cofiniteness was also relaxed slightly, allowing applications to some interesting examples such as affine VOAs at admissible level.

One of the features of \cite{meCMP} is the appearance of odd trace functions (see Section \ref{meat} for the definition) which are to be included alongside the more usual (super)trace functions in order to achieve modular invariance. Although similar in many ways, these odd traces differ from (super)traces in that they act nontrivially on odd elements of a vector superspace, whereas the (super)trace must always vanish on such elements. The results of \cite{meCMP} are reviewed in Section \ref{meat} (for simplicity in the special case of Ramond twisted modules).

In the present work, in Section \ref{ex4}, we compute an odd trace function for a particular example: the $N=1$ superconformal minimal model of central charge $c = -21/4$. We evaluate the odd trace function on the superconformal generator (which is an odd element of conformal weight $3/2$), using the strong constraint of its modular invariance. The odd trace function in question equals the weight $3/2$ modular form $\eta(\tau)^3$, where $\eta(\tau)$ is the well known Dedekind eta function.

We then give a different proof of this equality (up to an ambiguity of signs) using a BGG resolution and some simple combinatorics. The result is a representation theoretic interpretation of the classical identity 
\begin{align*}
\eta(\tau)^3 = q^{1/8} \sum_{n \in \Z} (4n+1) q^{n(2n+1)}
\end{align*}
similar in spirit, but a little different, to the celebrated proof coming from the affine Weyl-Kac denominator identity (\cite{IDLA} Chapter 12).

\emph{Acknowledgements:} I would like to warmly thank the organisers of the conference `Lie Superalgebras' at Universit\`{a} di Roma Sapienza where this work was presented. Also to express my gratitude to IMPA and to the IHES where the writing of this paper was completed.

\section{Definitions}\label{definitions}

For us a \underline{vertex superalgebra} \cite{Kac}, \cite{FBZ} is a quadruple $V, \vac, T, Y$ where $V$ is a vector superspace, $\vac \in V$ an even vector, $T : V \rightarrow V$ an even linear map, and $Y : V \otimes V \rightarrow V((z))$, denoted $u \otimes v \mapsto Y(u, z)v = \sum_{n \in \Z} u_{(n)}v z^{-n-1}$, is also even. These data are to satisfy the following axioms.
\begin{itemize}
\item The unit identities $Y(\vac, z) = I_V$ and $Y(u, z)\vac|_{z=0} = u$.

\item The translation invariance identity $Y(Tu, z) = \partial_z Y(u, z)$.

\item The Cousin property that the three expressions
\[
Y(u, z)Y(v, w) x \quad \quad  p(u, v) Y(v, w) Y(u, z) x, \quad \text{and} \quad Y(Y(u, z-w)v, w) x,
\]
which are elements of $V((z))((w))$, $V((w))((z))$, and $V((w))((z-w))$, are images of a single element of $V[[z, w]][z^{-1}, w^{-1}, (z-w)^{-1}]$ under natural inclusions into those three spaces.
\end{itemize}

An equivalent definition, more convenient for some applications, is the following. A vertex superalgebra is a triple $V, \vac, Y$ where these data are as above, but satisfy the following axioms.
\begin{itemize}
\item The unit identities $\vac_{(n)}u = \delta_{n, -1} u$, $u_{(-1)}\vac = u$ and $u_{(n)}\vac = 0$ for $n > 0$.

\item The Borcherds identity (also known, in a different notation, as the Jacobi identity)
\[
B(u, v, x; m, k, n) = 0 \quad \text{for all $u, v, x \in V$, $m, k, n \in \Z$},
\]
where
\begin{align*}
B(u, v, x; m, k, n) &= \sum_{j \in \Z_+} \binom{m}{j} (u_{(n+j)}v)_{(m+k-j)} x \\
&- \sum_{j \in \Z_+} (-1)^j \binom{n}{j} \left[ u_{(m+n-j)} v_{(k+j)} - (-1)^n p(u, v) v_{(k+n-j)} u_{(m+j)} \right] x.
\end{align*}
\end{itemize}
A \underline{vertex algebra} is a purely even vertex superalgebra.

Let $V$ be a vertex superalgebra. A $V$-module is a vector superspace $M$ with a vertex operation $Y^M : V \otimes M \rightarrow M((z))$ such that
\begin{align}\label{moduleborcherds}
Y^M(\vac, z) = I_M, \quad \text{and} \quad B(u, v, x; m, k, n) = 0
\end{align}
for all $u, v \in V$, $x \in M$, $m, k, n \in \Z$.

For the present theory we require the extra structure of a \underline{conformal vector}. This is a vector $\om \in V$ such that its modes $L_n = \om_{(n-1)} \in \en V$ furnish $V$ with a representation of the Virasoro algebra, i.e.,
\[
[L_m, L_n] = (m-n) L_{m+n} + \delta_{m, -n} \frac{m^3-m}{12} c
\]
(here $c \in \C$ is an invariant of $V$ called the central charge), $L_0$ is diagonal with real eigenvalues bounded below, and $L_{-1} = T$. We call a vertex superalgebra with conformal vector a \underline{vertex operator superalgebra} or VOSA, and we use the term VOA to distinguish the purely even case.

A $V$-module $M$ is called a \underline{positive energy} module if $L_0 \in \en M$ acts diagonally with eigenvalues bounded below. In particular $V$ is a positive energy $V$-module. The eigenvalues of $L_0 \in \en V$ are called \underline{conformal weights}, and if $L_0 u = \D u$ we write
\[
Y(u, z) = \sum_{n \in \Z} u_{(n)} z^{-n-1} = \sum_{n \in -\D + \Z} u_n z^{-n-\D}
\]
(so that $u_n = u_{(n-\D+1)}$). The \underline{zero mode} $u_0 \in \en M$ attached to $u \in V$ is special because it commutes with $L_0$ and thus preserves the eigenspaces of the latter. A VOSA $V$ is said to be \underline{rational} if its category of positive energy modules is semisimple, i.e., it contains finitely many irreducible objects, and any object is isomorphic to a direct sum of irreducible ones.

The condition of $C_2$-cofiniteness is an important finiteness condition of vertex (super)algebras introduced by Zhu. We say that $V$ is \underline{$C_2$-cofinite} if
\[
\dim \left( V / V_{(-2)}V \right) < \infty.
\]

\section{The Theorem of Zhu}\label{zhusection}

Now we come to the theorem of Zhu \cite{Zhu}.
\begin{thm}[Zhu] \label{zhuthm}
Let $V$ be a VOA (i.e., purely even VOSA) such that:
\begin{itemize}
\item $V$ is rational,

\item the conformal weights of $V$ lie in $\Z_+$,

\item $V$ is $C_2$-cofinite.
\end{itemize}
We associate to each irreducible positive energy module $M$, and $u \in V$, the series
\[
S_M(u, \tau) = \tr_M u_0 q^{L_0 - c/24},
\]
convergent for $q = e^{2\pi i \tau}$ of modulus less than $1$. There is a grading $V = \oplus_{\nabla \in \Z_+} V_{[\nabla]}$ such that for $u \in V_{[\nabla]}$ the span $\CC(u)$ of the (finitely many) functions $S_M(u, \tau)$ defined above is modular invariant of weight $\nabla$, i.e.,
\[
(c\tau + d)^\nabla f\left(\frac{a\tau+b}{c\tau+d}\right) \in \CC(u) \quad \text{for all $f(\tau) \in \CC(u)$ and $\slmat \in SL_2(\Z)$}.
\]
\end{thm}

Here is an outline of the proof of Zhu's theorem.
\begin{enumerate}
\item Introduce a space $\CC$ of maps $S(u, \tau) : V \times \HH \rightarrow \C$ linear in $V$ and holomorphic in $\HH = \{\tau \in \C | \text{Im}{\tau} > 0\}$ satisfying certain axioms, this space is called the `conformal block' of $V$ or the space of conformal blocks of $V$. The definition of $\CC$ can be understood in terms of elliptic curves and their moduli \cite{FBZ}.

\item It is automatic from its definition that $\CC$ admits an action of the group $SL_2(\Z)$, namely,
\[
[S \cdot A](u, \tau) = (c\tau + d)^{-\nabla_u} S(u, A\tau),
\]
where $\nabla$ is the grading mentioned above.

\item It is proved by direct calculation that $\tr_M u_0^M q^{L_0 - c/24}$ is a conformal block (at least as a formal power series).

\item Using the $C_2$-cofiniteness condition, one shows that any fixed $S \in \CC$ satisfies some differential equation, and consequently is expressible as a power series in $q$ (whose coefficients are linear maps $V \rightarrow \C$).

\item \label{zhustep} The lowest order coefficient $C_0 : V \rightarrow \C$ in the series expansion factors to a certain quotient $\zhu(V)$ of $V$. This quotient has the structure of a unital associative algebra, and $C_0$ is symmetric, i.e., $C_0(ab) = C_0(ba)$.

\item There is a natural bijection
\[
\text{irreducible positive energy $V$-modules} \longleftrightarrow \text{irreducible $\zhu(V)$-modules},
\]
and so if $V$ is rational, $\zhu(V)$ is finite dimensional semisimple. Thus we can write $C_0 = \sum_N \al_N \tr_N$ where the sum is over irreducible $\zhu(V)$-modules and $\al_N \in \C$.

\item Write the corresponding sum $\sum_N \al_N S_M$ where $M$ is the $V$-module associated to $N$. Subtract this conformal block from $S$.

\item One can repeat the process and show that $S$ is exhausted by trace functions in a finite number of steps.
\end{enumerate}

\section{Generalisation to the Supersymmetric Case, and to Rational Conformal Weights}\label{meat}

Results described in this section are drawn from \cite{meCMP}.

Many examples of interest, especially in the supersymmetric case, are graded by noninteger conformal weights. So it is first necessary (and of independent interest) to relax the condition of integer conformal weights in Theorem \ref{zhuthm}. Therefore let $V$ be a VOA whose conformal weights lie in $\Q$ (and are bounded below) rather than in $\Z_+$, and which satisfies the other conditions of the theorem. Then the trace functions $S_M(u, \tau)$ and their span $\CC(u)$ are defined as before. There exists a certain \emph{rational} grading $V = \oplus_{\nabla \in \Q} V_{[\nabla]}$ in place of the usual integer grading. It is then true that for $u \in V_{[\nabla]}$ the space $\CC(u)$ is invariant under the weight $\nabla$ action\footnote{The precise definition of the action involves choices of roots of unity in general. See \cite{meCMP} for details.}, not of $SL_2(\Z)$, but of its congruence subgroup
\begin{align*}
\G_1(N) = \{\slmat \in SL_2(\Z) | b \equiv 0 \bmod N, \,\text{and}\,\, a \equiv d \equiv 1 \bmod N\}.
\end{align*}
Here $N$ is the least common multiple of the denominators of conformal weights of vectors in $V$. This number is finite because of the condition of $C_2$-cofinitness.

It is possible to achieve invariance under the whole of $SL_2(\Z)$ by altering our definition of $V$-module. Define a \underline{Ramond twisted} $V$-module to be a vector superspace $M$ together with fields
\[
Y^M(u, z) = \sum_{n \in -\D_u + \Z} u_{(n)} z^{-n-1} = \sum_{n \in \Z} u_n z^{-n-\D}
\]
satisfying (\ref{moduleborcherds}) for all $m \in -\D_u + \Z$, $k \in -\D_v + \Z$, $n \in \Z$. Notice that the ranges of indices of the modes are modified so that $u \in V$ always possesses integrally graded modes $u_n \in \en M$, and in particular always possesses a zero mode. Let us call a VOA \underline{Ramond rational} if its category of positive energy Ramond twisted modules is semisimple.

Let $V$ be a Ramond rational, $C_2$-cofinite VOA with rational conformal weights bounded below. Attach the trace function
\[
S_M(u, \tau) = \tr_M u_0 q^{L_0 - c/24}
\]
to $u \in V$ and $M$ an irreducible positive energy Ramond twisted $V$-module, and let $\CC(u)$ be the span of all such trace functions. Then for $u \in V_{[\nabla]}$ the space $\CC(u)$ is invariant under the weight $\nabla$ action of the full modular group $SL_2(\Z)$.

The previous paragraph stated a result for VOAs. Upon passage from VOAs to VOSAs, one might expect the claim to hold with trace functions simply replaced by supertrace functions $S_M(u, \tau) = \str_M u_0 q^{L_0 - c/24}$. This would be true, except for an interesting subtlety which can be traced to Step \ref{zhustep} of the proof outline given in Section \ref{zhusection}.

In the present situation it is appropriate to replace the usual Zhu algebra with a certain superalgebra (which we also refer to as the Zhu algebra and denote $\zhu(V)$) introduced in the necessary level of generality in \cite{DK}. If $V$ is Ramond rational then $\zhu(V)$ is finite dimensional and semisimple. The lowest coefficient $C_0$ of the series expansion of a conformal block now descends to a supersymmetric function on $\zhu(V)$, i.e., $C_0(ab) = p(a, b)C_0(ba)$.

The classification of pairs $(A, \vp)$, where $A$ is a finite dimensional simple superalgebra (over $\C$) and $\vp$ is a supersymmetric function on $A$, is as follows.
\begin{itemize}
\item $A = \en(N)$ for some finite dimensional vector superspace $N$, and $\vp$ is a scalar multiple of $\str_N$.

\item $A = \en(P)[\theta] / (\theta^2 - 1)$ where $P$ is a vector space and $\theta$ is an odd indeterminate, and $\vp$ is a scalar multiple of $a \mapsto \tr_P(a \theta)$.
\end{itemize}
The first case is the analogue of the usual Wedderburn theorem. The superalgebra of the second case is known as the \underline{queer superalgebra} and is often denoted $Q_n$ (where $n = \dim P$). Clearly we have
\[
Q_n \cong \left\{ \twobytwo X Y Y X | X, Y \in \text{Mat}_n(\C) \right\},
\]
and the unique up to a scalar factor supersymmetric function on $Q_n$ is $\left(\begin{smallmatrix}{X}&{Y}\\{Y}&{X}\\ \end{smallmatrix}\right) \mapsto \tr Y$, which is known as the \underline{odd trace}.

Roughly speaking modular invariance will hold for $\CC(u)$ defined as the span of supertrace functions together with apropriate analogues of odd trace functions. More precisely:
\begin{defn}\label{queertracedefn}
Let $V$ be a Ramond rational, $C_2$-cofinite VOSA with rational conformal weights bounded below. Let $A$ be a simple component of $\zhu(V)$, $N$ the corresponding unique $\Z_2$-graded irreducible module, and $M$ the corresponding irreducible positive energy Ramond twisted $V$-module. If $A \cong \en(P)[\theta]/(\theta^2-1)$ is queer then let $\Theta : M \rightarrow M$ denote the lift to $M$ of the map $\theta : N \rightarrow N$ of multiplication by $\theta$.
 In this case we define the \underline{odd trace function}
\[
S_M(u, \tau) = \tr_M u_0 \Theta q^{L_0 - c/24}.
\]
If $A$ is not queer then we define the \underline{supertrace function}
\[
S_M(u, \tau) = \str_M u_0 q^{L_0 - c/24}.
\]
\end{defn}
Now we can state the main theorem: it is Theorem 1.3 of \cite{meCMP} applied to the special case of untwisted characters of Ramond twisted $V$-modules.
\begin{thm}[\cite{meCMP}]\label{mythm}
Let $V$ be a VOSA as in Definition \ref{queertracedefn}, and let $\CC(u)$ be the span of the supertrace functions and odd trace functions attached to all irreducible positive energy Ramond twisted $V$-modules. There exists a grading $V = \oplus_{\nabla \in \Q} V_{[\nabla]}$ such that for $u \in V_{[\nabla]}$ the space $\CC(u)$ is invariant under the weight $\nabla$ action of $SL_2(\Z)$.
\end{thm}
In the next two sections we view some examples of odd trace functions.

\section{Example: The Neutral Free Fermion}

This example is considered more fully in \cite{meCMP}, the interested reader may refer there for further details.

We consider the Lie superalgebras $A^\text{tw}$ (resp. $A^\text{untw}$)
\[
(\oplus_{n} \C \psi_n) \oplus \C 1
\]
where the direct sum ranges over $n \in 1/2 + \Z$ (resp. $n \in \Z$). Here the vector $1$ is even, $\psi_n$ is odd. The commutation relations in both cases are
\[
[\psi_{m}, \psi_{n}] = \delta_{m, -n} 1.
\]

We introduce the Fock module
\[
V = U(A^\text{tw}) \otimes_{U(A^\text{tw}_+)} \C \vac,
\]
where
\[
A^\text{tw}_+ = \C 1 \oplus (\oplus_{n \geq 1/2} \C \psi_n)
\]
and $\C \vac$ is the $A^\text{tw}_+$-module on which $1$ acts as the identity and $\psi_n$ acts trivially.

It is well known \cite{Kac} that $V$ can be given the structure of a VOSA. The Virasoro element is $\om = \frac{1}{2} \psi_{-3/2}\psi_{-1/2}\vac$ with central charge $c = 1/2$. The vector $\psi = \psi_{-1/2}\vac$ has conformal weight $1/2$ and associated vertex operator
\[
Y(\psi, z) = \sum_{n \in 1/2 + \Z} \psi_n z^{-n-1/2}.
\]

Ramond twisted $V$-modules are, in particular, modules over the untwisted Lie superalgebra $A^\text{untw}$. In fact the unique irreducible positive energy Ramond twisted $V$-module is 
\[
M = U(A^\text{untw}) \otimes_{U(A^\text{untw}_+)} (\C v + \C \ov{v})
\]
where
\[
A^\text{untw}_+ = \C 1 \oplus (\oplus_{n \geq 0} \C \psi_n)
\]
and $\C v + \C \ov{v}$ is the $A^\text{untw}_+$-module on which $1$ acts as the identity, $\psi_n$ acts trivially for $n > 0$, $\psi_0 v = \ov{v}$ and $\psi_0 \ov{v} = v/2$. We note that $V$ is $C_2$-cofinite and Ramond rational (as well as being rational).

The (Ramond) Zhu algebra of $V$ is explicitly isomorphic to the queer superalgebra $Q_1 = \C[\theta] / (\theta^2 = 1)$ via the map $[\vac] \mapsto 1$, $[\psi] \mapsto \sqrt{2} \theta$. Thus the lowest graded piece $M_0 = \C v + \C \ov{v}$ of $M$ is a $Q_1$-module.

The corresponding odd trace function is
\[
S_M(u, \tau) = \tr_M u_0 \Theta q^{L_0-c/24}
\]
where $\Theta : M \rightarrow M$ is as in Definition \ref{queertracedefn}. By unwinding that definition we see that $\Theta : \mathbf{m} w \mapsto \mathbf{m} \psi_0 w$, where $w$ is $v$ or $\ov{v}$, and the monomial $\mathbf{m} \in U(A^\text{untw} / A^\text{untw}_+)$.

The odd trace function $S_M(u, \tau)$ vanishes on $u = \vac$ (indeed on all even vectors of $V$), but it acts nontrivially on the odd vector $\psi$ (which is pure of Zhu weight $1/2$). Therefore $S_M(\psi, \tau)$ must be a modular form on $SL_2(\Z)$ of weight $1/2$ (with possible multiplier system).

Indeed one may verify that $\psi_0 \Theta$ acts as $(-1)^{\text{length}(\mathbf{m})}$ on the monomial vector $\mathbf{m} w$, and so
\begin{align*}
S_M(\psi, \tau) &= q^{-c/24} q^{L_0|{M_0}} (1 - q^1) (1 - q^2) \cdots \\
&= q^{1/24} \prod_{n=1}^\infty (1-q^n) = \eta(\tau)
\end{align*}
(here we have used that $c = 1/2$, and that $L_0|_{M_0} = 1/16$). We have recovered the well known Dedekind eta function $\eta(\tau)$ which is indeed a modular form on $SL_2(\Z)$ of weight $1/2$.

\section{Example: The $N=1$ Superconformal Algebra} \label{ex4}

First we recall the definition of the Neveu-Schwarz Lie superalgebra $\nstw$, and its Ramond-twisted variant $\nsuntw$ (which is often called the Ramond superalgebra).
\begin{defn}
As vector superspaces the Lie superalgebras $\nstw$ (resp. $\nsuntw$) are
\[
(\oplus_{n \in \Z} \C L_n) \oplus (\oplus_{m} \C G_n) \oplus \C C
\]
where the direct sum ranges over $m \in \Z$ (resp. $m \in 1/2 + \Z$). Here $C$ and $L_n$ are even, $G_m$ is odd. The commutation relations in both cases are
\begin{align} \label{NS}
\begin{split}
[L_m, L_n] &= (m-n)L_{m+n} + \frac{m^3-m}{12} \delta_{m, -n} C, \\
[G_m, L_n] &= (m - \frac{n}{2}) G_{m+n}, \\
[G_m, G_n] &= 2L_{m+n} + \frac{1}{3}(m^2 - \frac{1}{4}) \delta_{m, -n} C,
\end{split}
\end{align}
with $C$ central.
\end{defn}
As usual we introduce the Verma $\nstw$-module
\[
M^{\text{tw}}(c, h) = U(\nstw) \otimes_{U(\nstw_+)} \C v_{c, h}
\]
where 
\begin{align*}
\nstw_+ = \C C + \oplus_{n \geq 0} \C L_n + \oplus_{m \geq 1/2} \C G_m,
\end{align*}
and $\C v_{c, h}$ is the $\nstw_+$-module on which $C$ acts by $c \in \C$, $L_0$ acts by $h \in \C$, and higher modes act trivially. It is well known \cite{Kac}, \cite{KacWang} that the quotient $\NS^c = M^\text{tw}(c, 0) / U(\nstw) G_{-1/2} v_{c, 0}$ is a VOSA of central change $c$, as is the irreducible quotient $\NS_c$.

We shall also require the (generalised) Verma $\nsuntw$-modules
\[
M(c, h) = U(\nsuntw) \otimes_{U(\nsuntw_+)} S_{c, h}
\]
and their irreducible quotients $L(c, h)$ (we omit the superscript $\text{untw}$) where
\begin{align*}
\nsuntw_+ = \C C + \oplus_{n \geq 0} \C L_n + \oplus_{m \geq 0} \C G_m,
\end{align*}
and $S_{c, h}$ is the $\nsuntw_+$-module characterised by:
\begin{align*}
\begin{array}{ll}
\text{$S_{c, h} = \C v_{c, h}$ with $G_0 v_{c, h} = 0$} & \text{if $h = c/24$}, \\
\text{$S_{c, h} = \C v_{c, h} + \C G_0 v_{c, h}$} & \text{if $h \neq c/24$}, \\
\end{array}
\end{align*}
with $C = c$, $L_0 = h$, and positive modes acting trivially in both cases.

A clear summary of the representations of $\NS^c$ and $\NS_c$ can be found in \cite{supermilas}. Here we focus on the Ramond twisted representations and shall often omit the adjective `Ramond twisted'. Generically $\NS_c = \NS^c$ is irreducible and all the $\nsuntw$-modules $L(c, h)$ acquire the structure of positive energy $\NS_c$-modules. For certain values of $c$ though, $\NS_c$ is a nontrivial quotient of $\NS^c$ and the irreducible positive energy $\NS_c$-modules are finite in number and are all of the form $L(c, h)$. In fact, $\NS_c$ is a (Ramond) rational VOSA when
\begin{align*}
c = c_{p, p'} = \frac{3}{2}\left(1 - \frac{2(p'-p)^2}{pp'}\right)
\end{align*}
for $p, p' \in \Z_{>0}$ with $p < p'$, $p' - p \in 2\Z$ and $\gcd(\frac{p'-p}{2}, p) = 1$. In this case the irreducible positive energy $\NS_c$-modules are precisely the $\nsuntw$-modules $L(c, h)$ where
\begin{align*}
h = h_{r, s} = \frac{(rp'-sp)^2 - (p'-p)^2}{8pp'} + \frac{1}{16}
\end{align*}
for $1 \leq r \leq p-1$ and $1 \leq s \leq p'-1$ with $r-s$ odd.

Let $c$ be one of these special values from now on. The irreducible $\zhu(\NS_c)$-modules are precisely the lowest graded pieces of the modules $L(c, h)$ introduced above. The lowest graded piece is of dimension $1$ if $h = c/24$ (there is clearly at most one such module for any fixed value of $c$), and is of dimension $1|1$ if $h \neq c/24$. It is known that $\zhu(\NS_c)$ is supercommutative (it is a quotient of $\zhu(\NS^c) \cong \C[x, \theta] / (\theta^2 - x + c/24)$ where $x$ is even and $\theta$ odd). Therefore the simple components of $\zhu(\NS_c)$ with the $1|1$-dimensional modules are all copies of the queer superalgebra $Q_1$, and the component with the $1$-dimensional module (if it exists) is $\C$.

Let us consider the case $c = -21/4$ (so $p=2$, $p'=8$) for which the two irreducible positive energy modules are $M_i = L(c, h_i)$ where $h_1 = -3/32$ and $h_2 = -7/32 = c/24$. The first of these is the unique queer module. Theorem \ref{mythm} tells us that the supertrace function $S_{M_2}(u, \tau)$ and the odd trace function $S_{M_1}(u, \tau)$ together span an $SL_2(\Z)$-invariant space whose weight is the Zhu weight of $u$. Assume further that $u \in V$ is odd. Then $S_{M_2}(u, \tau)$ vanishes as the supertrace of an odd element. But $S_{M_1}(u, \tau)$ need not vanish, and it will be a modular form (with multiplier system).

Unwinding Definition \ref{queertracedefn} we see that $\Theta : \mathbf{m} v_{c, h} \mapsto \mathbf{m} G_0 v_{c, h}$ where $\mathbf{m} \in U(\nsuntw / \nsuntw_+)$, and
\[
S_{M_1}(u, \tau) = \tr_{M_1} u_0 \Theta q^{L_0 - c/24}.
\]

The VOSA $\NS_c$ possesses a distinguished element $\nu = G_{-3/2}\vac$ of conformal weight $3/2$, it satisfies $\nu_0 = G_0$. It turns out that $\nu$ is of pure Zhu weight $3/2$ and so by the above remarks
\[
F(\tau) := \tr_{M_1} G_0 \Theta q^{L_0 - c/24}
\]
is a modular form of weight $3/2$. On the top level of $M_1$, $G_0 \Theta = G_0^2 = h - c/24 = 1/8$, so the top level contribution to $F(\tau)$ is $\frac{1}{4}q^{1/8}$. This is already enough information to determine $F(\tau)$ completely. The cube of the Dedekind eta function is $q^{1/8}$ times an ordinary power series in $q$, so the quotient $f(\tau) = F(\tau) / \eta(\tau)^3$ is a holomorphic modular form of weight $0$ for $SL_2(\Z)$, possibly with a multiplier system. Since the $q$-series of $f$ has integer powers of $q$ we have $f(T\tau) = f(\tau)$, and since $S^2 = 1$ we have only the possibilities $f(S\tau) = \pm f(\tau)$. But in $SL_2(\Z)$ we have the relation $(ST)^3 = 1$, so if $S$ acted by $-1$ on $f$ we would have $f(\tau) = f(-T^3\tau) = -f(\tau)$. Hence $f(S\tau) = f(\tau)$ and, since it is a genuine holomorphic modular form on $SL_2(\Z)$, we have $f(\tau) = 1$. Thus
\begin{align} \label{calcofS}
F(\tau) = \eta(\tau)^3 / 4.
\end{align}

We next compute $F(\tau)$ using representation theory. We obtain (up to some undetermined signs) the following well known classical identity of Jacobi
\begin{align}\label{jaceta}
\eta(\tau)^3 = q^{1/8} \sum_{n \in \Z} (4n+1) q^{n(2n+1)}.
\end{align}

We begin by considering the trace of $G_0 \Theta q^{L_0}$ on the Verma module $M(c, h)$. The action of $G_0 \Theta$ on the monomial
\begin{align*}
\mathbf{m} v = L_{m_1} \cdots L_{m_s} G_{n_1} \cdots G_{n_t} v,
\end{align*}
where $m_1 \leq \ldots \leq m_s \leq -1$, $n_1 < \ldots < n_t \leq -1$, and $v$ is $v_{c, h}$ or $G_0 v_{c, h}$, looks like
\[
\mathbf{m} G_0 v \mapsto G_0 \mathbf{m} G_0 v = (-1)^t \mathbf{m} G_0^2 v + \text{reduced terms}.
\]
Reduced terms resulting from a single use of the commutation relations are of the same length as $\mathbf{m}$, but contain different numbers of the symbols $L$ and $G$. Reduced terms resulting from more than one use of the commutation relations are strictly shorter than $\mathbf{m}$. Therefore none of these terms contribute to the trace. Consider monomials $\mathbf{m}$ as above with a fixed value of $N = \sum_{i=1}^s m_i + \sum_{j=1}^t n_j$. A simple generating function argument shows that if $N > 0$ then the number of such monomials with $t$ even is the same as the number with $t$ odd. Thus the only nonzero term in $\tr G_0 \Theta q^{L_0}$ is the leading term.

It is known that $L(c = -21/4, h_0 = -3/32)$ has a BGG resolution
\[
0 \leftarrow L(c, h_0) \leftarrow M(c, h_0) \leftarrow M(c, h_1) \oplus M(c, h_{-1}) \leftarrow M(c, h_2) \oplus M(c, h_{-2}) \leftarrow \cdots,
\]
where $h_n = -3/32 - n(2n+1)$ for all $n \in \Z$, and that all $M(c, h_k)$ are naturally embedded in $M(c, h_0)$ \cite{IoharaKoga}. We therefore identify each Verma module with its image in $M(c, h_0)$. The trace we seek is given as an alternating sum over the terms in the resolution.

From this we see already that the only nonzero coefficients in the $q$-expansion of $\eta(\tau)^3$ must be for powers $q^{1/8 - n(2n+1)}$. We can also easily determine the coefficients up to a sign. Indeed we have $\Theta^2 = h_0 - c/24 = 1/8$, while the operator $G_0|_{M(c, h_k)}$ preserves the top piece $S_k$ of $M(c, h_k)$ and squares to $h_k^2 - c/24$. Therefore the operator $(G_0 \Theta)|_{S_k}$ (which is diagonal on the $1|1$-dimensional space $S_k$) squares to
\[
(h_k-c/24)^2/8 = \left[(4k+1)/8\right]^2.
\]
This matches perfectly with (\ref{jaceta}). To determine the signs of the coefficients directly it seems to be necessary to know some further information about the singular vectors, it would be nice to find a simpler derivation.

Of course similar arguments may be applied to other rational $\NS_c$ and their modules. If $L_0$ happens to take the value $c/24$ on one of the levels of a module then the arguments potentially become more intricate.

\end{document}